\documentclass[letterpaper, preprint, paper,11pt]{AAS}	

\usepackage{bm}
\usepackage{amsmath}
\usepackage{subcaption}
\usepackage{booktabs}
\usepackage[colorlinks=true, pdfstartview=FitV, linkcolor=black, citecolor= black, urlcolor= black]{hyperref}
\usepackage{overcite}
\usepackage{footnpag}			      	
\usepackage{amssymb}
\usepackage{caption}
\usepackage[ruled, linesnumbered]{algorithm2e}
\usepackage{graphicx}
\usepackage{float}

\PaperNumber{XX-XXX}
\usepackage[super,comma,sort&compress]{natbib}

\makeatletter
\newcommand*{\citenst}[2][]{%
  \begingroup
  \let\NAT@mbox=\mbox
  \let\@cite\NAT@citenum
  \let\NAT@space\NAT@spacechar
  \let\NAT@super@kern\relax
  \renewcommand\NAT@open{[}%
  \renewcommand\NAT@close{]}%
  \citet[#1]{#2}%
  \endgroup
}
\makeatother

\begin{document}

\title{Nonlinear Reachable Set Computation and Model Predictive Control for Safe Hypersonic Re-entry of Atmospheric Vehicles}

\author{Jinaykumar Patel\thanks{Graduate Student, Mechanical and Aerospace Engineering, The University of Texas at Arlington, Arlington, TX 76019 } and Kamesh Subbarao\thanks{Professor, Mechanical and Aerospace Engineering, The University of Texas at Arlington, Arlington, TX 76019}
}

\maketitle{} 		

\begin{abstract}
This paper investigates the application of reachability analysis to the re-entry problem faced by vehicles entering Earth’s atmosphere. The study delves into the time evolution of reachable sets for the system, particularly when subject to nonlinear implicit controls, given the potential damage from the intense heat generated during hypersonic re-entry. Our proposed methodology leverages zonotopes and constrained zonotopes to ensure compliance with safety specifications. Furthermore, we utilize Model Predictive Control for detailed trajectory planning. To substantiate our methodology, we provide detailed simulations that not only tackle nonlinear re-entry scenarios but also illustrate trajectory planning using MPC.
\end{abstract}

\section{Introduction}
Recent advancements in human spaceflight programs have amplified interest in asteroid mining, space tourism, and the establishment of outposts on celestial entities~\cite{arslantacs2016safe}. These endeavors pose intricate challenges, such as the re-entry, descent, and landing of spacecraft~\cite{rajnishReentry,katiReentry}. The hypersonic re-entry phase of flight poses significant challenges for vehicles due to the intense heat generated, which can cause damage to the vehicle's structure and materials if not properly managed. Precise guidance and control are essential for the mission's success, mainly since hypersonic re-entry includes dynamics with significant nonlinearities, state and control constraints, and parameter uncertainties. Additionally, a safety system is required to ensure the system's robustness. A map of attainable states of the dynamical system is necessary to mitigate concerns related to re-planning trajectories to alternate landing sites. This is where reachability analysis comes in, as it enables us to obtain a comprehensive set of all possible states for the system, given its initial condition and admissible control inputs.

Reachable sets are the set of all states that a system can attain given the admissible inputs, and they have been studied widely in literature due to their numerous applications. These include collision avoidance~\cite{zhou2015reachable}, obstacle avoidance~\cite{malone2017hybrid}, and goal satisfaction, where computing the reachable set is used to confirm safety. Motion planning~\cite{mcmahon2014sampling} - utilizes computation of optimal reachable set for effective planning. In the stochastic realm, reachability analysis is typically carried out, assuming the system dynamics are linear. However, unlike linear systems, it is only possible to calculate exact reachable sets for a nonlinear system if the system is integrable and free of uncertainties. Therefore, numerical approximation methods have been proposed to compute the reachable sets~\cite{lofberg_yalmip_2004, girard2005reachability, althoff_implementation_2016, kurzhanskiy_ellipsoidal_2006}. These numerical methods generally involve the use of ellipsoidal sets and zonotopes~\cite{althoff_reachability_2008, filippova_ellipsoidal_2017}.  

Spacecraft re-entry into Earth's atmosphere is fraught with challenges, most notably the intense heating caused by friction between the spacecraft's surface and atmospheric particles. It's of utmost importance to manage this heating effectively, as failure to do so could lead to structural damage or even the complete destruction of the spacecraft. Model Predictive Control (MPC) offers a promising solution to address this challenge. After proving its mettle in diverse applications, from the process industry to automotive sectors, MPC has been progressively recognized as aptly suited for the intricate challenges of space missions. Its effectiveness and adaptability have been a subject of rigorous studies, as seen in references~\cite{schurmann2018reachset, 7330732}. MPC's methodology hinges on its ability to forecast future spacecraft outputs and control actions based on past data. This approach ensures that trajectories are optimized while always acknowledging the inherent constraints. Moreover, a standout feature of MPC is its resilience in compensating for modeling uncertainties. The inherent feedback control within MPC bolsters its precision, making it an invaluable tool in such critical scenarios. Its capacity to manage specific system constraints is explicit, positioning it as a preferred choice for providing guidance and control, especially when integrated with relative motion dynamic models for predictive purposes.

In the specific context of spacecraft re-entry, MPC ensures that the vehicle adheres to a reference trajectory, constrained within a predetermined reachable set tube. If any uncertainties trigger a deviation from this trajectory, the control mechanisms within MPC are primed to guide the spacecraft back into the safety of the tube. This becomes particularly crucial when the spacecraft operates near its heat load limits. Not only does the control mechanism prevent the spacecraft from entering perilous zones, but it also aids in course correction if the spacecraft has inadvertently exceeded these limits. The reliability of MPC is underscored by the assurance that, once inside the constrained zonotope, the spacecraft remains within the safety parameters, ensuring any reference trajectory is maintained while keeping the heating rate within permissible bounds.

This paper builds upon the previous work in [~\citenum{jinay_JGCD}] by extending the analysis to nonlinear systems and focusing specifically on the re-entry problem. While the previous work primarily dealt with linear and perturbed systems, we delve deeper into the intricacies of nonlinear systems and their reachability analysis. To achieve this, we present the computation methods for two sets of types: zonotopes and constrained zonotopes, with the mathematical formulation of these sets being based on the earlier work. Our analysis also includes a detailed examination of the impact of heating rate on the reachable sets of the re-entry vehicle, with constrained zonotopes ensuring that safety specifications are met. Overall, this paper provides a comprehensive reachability analysis in the context of the re-entry problem, building upon existing work and presenting new insights into the complexities of nonlinear systems.
\section{Problem Statement}
\label{sec:prob}
Before computing the reachable sets, we first review the model equations. Various versions and problem formulations of atmospheric re-entry exist in the literature. In this context, we consider the nonlinear equations of motion for a hypersonic vehicle during atmospheric re-entry~\cite{graichen2008constructive, jain2020computationally}:
\begin{equation} \label{eq:16}
\begin{aligned}
\dot{h} &=v \sin \gamma \\
\dot{v} &=\frac{-D(h, v, \alpha)}{m}-g(h) \sin \gamma \\ 
\dot{\gamma} &=\frac{L(h, v, \alpha)}{m v} \cos \beta+\cos \gamma\left(\frac{v}{R_{e}+h}-\frac{g(h)}{v}\right) \\
\dot{\theta} &=\frac{v}{R_{e}+h} \cos \gamma \sin \psi \\
\dot{\psi} &=\frac{L(h, v, \alpha)}{m v} \frac{\sin \beta}{\cos \gamma}-\frac{v}{R_{e}+h} \cos \gamma \cos \psi \tan \theta \\
\dot{\phi} &=\frac{v}{R_{e}+h} \frac{\cos \gamma \cos \psi}{\cos \theta}
\end{aligned}
\end{equation}
where $h, v$ and $\gamma$ represent the altitude, velocity, and flight-path angle of the vehicle. $\theta, \phi$ and $\psi$ is latitude, longitude, and heading angle, respectively. The state vector $\mathbf{x}$ is defined as follows:
\begin{equation} \label{eq:17}
\mathbf{x}=[h~ v~ \gamma~ \theta~ \psi~ \phi]^{T}
\end{equation}
We also consider the aerodynamic and atmospheric forces acting on the vehicle during re-entry. The atmospheric density profile is typically modeled by the exponential barometric formula:
\begin{equation} \label{density}
	\rho(h)=\rho_{0} \exp \left[-h / h_{r}\right]
\end{equation}
where $\rho_{0}$ is the atmospheric density at sea level, and $h_r$ is the atmospheric scale height. The gravitational force on the vehicle depends on altitude and is given by:
\begin{equation}
	g(h)=G/\left(R_{e}+h\right)^{2}
\end{equation}
where $G$ is the universal gravitational constant and $R_e$ is the radius of the Earth/Planet.

The lift and drag forces acting on the vehicle can be modeled as functions of the air density, dynamic pressure, reference area, lift and drag coefficients, and vehicle orientation. We use the following models \cite{jain2020computationally}:
\begin{equation} \label{eq:18}
	\begin{array}{l}
		L(h, v, \alpha)=\frac{1}{2} C_{L}(\alpha) S \rho(h) v^{2}, \quad C_{L}(\alpha)=a_{0}+a_{1} \alpha \\[10pt]
		D(h, v, \alpha)=\frac{1}{2} C_{D}(\alpha) S \rho(h) v^{2}, \quad C_{D}(\alpha)=b_{0}+b_{1} \alpha+b_{2} \alpha^{2}
	\end{array}
\end{equation}
where, $S$ is the reference area. $a_i$ and $b_i$ are constants used to calculate the lift and drag acting on the vehicle.

\subsection{Representation of the System and Control Parameters}

For the system articulated by Eq.~(\ref{eq:16}), we identify two essential control parameters: the angle of attack ($\alpha$) and the bank angle ($\beta$). Both angles exert significant influence on the aerodynamic characteristics of a hypersonic re-entry vehicle. Specifically, the angle of attack is the angular separation between the vehicle's velocity vector and its longitudinal axis. Meanwhile, the bank angle delineates the angle between the longitudinal axis of the vehicle and a horizontal plane.

The interplay of these angles can drastically influence the aerodynamic forces and torques on the re-entry vehicle. As an illustration, an elevated angle of attack may escalate aerodynamic drag, potentially altering the vehicle's trajectory and stability. In a similar vein, pronounced bank angles can induce potent rolling torques, which might further perturb the vehicle's stability and course. Nevertheless, adept manipulation of both angles can enhance attributes like lift and stability, equipping the vehicle for a more regulated and secure re-entry.

Within Eq.~(\ref{eq:16}), these control parameters manifest implicitly and in a nonlinear manner. The relationship can be succinctly encapsulated as:

\begin{equation}
	\label{eq:16compact}
	\dot{\mathbf{x}} = \mathbf{g}\left( \mathbf{x}, \mathbf{u}\right),
\end{equation}

where $\mathbf{x} \in \mathbb{R}^{n}$, $\mathbf{u} \in \mathbb{R}^{m}$, and $\mathbf{g}: \mathbb{R}^{n} \times \mathbb{R}^{m} \rightarrow \mathbb{R}^{n}$. In the context of Eq.~(\ref{eq:16}), dimensions are defined by $n=6$ and $m=2$, hence, $\mathbf{g}: \mathbb{R}^{6} \times \mathbb{R}^{2} \rightarrow \mathbb{R}^{6}$.

Considering the nonlinear dynamics presented in Eq.~\ref{eq:16compact}, we'll determine reachable sets, denoted $\mathcal{R}(t)$, over the time interval $t \in \left[t_0~ t_f\right]$, with $t_0$ and $t_f$ marking the commencement and termination of the interval. We adopt the definition of a reachable set $\mathcal{R}(\cdot)$ as postulated in [\citenum{kailath}].

\subsection{Reachability}
For time-invariant systems, a state $\mathbf{x}_1$ is deemed reachable if there exists an input capable of transitioning the system state from $\mathbf{x}(t_0)$ to $\mathbf{x}_1$ within a finite duration, denoted by $T = t_f - t_0$. Given a predetermined state $\mathbf{x}_1$, one can compute a distinct control history to effectuate this transition, provided that the system is inherently reachable. On the other hand, by employing the system's initial state $\mathbf{x}(t_0)$, the governing dynamics as portrayed in Eq.~\ref{eq:16compact}, and the bounded control vector $\mathbf{u}$, one can derive a set encompassing all feasible values of $\mathbf{x}(t)$ at any specific time $t$. This is achieved by methodically sampling the control vector over the interval $[t_0, t]$. Such a set is conventionally referred to as the \textit{Reachable Set}.

\section{Set Representations: A Primer}
In the context of geometric analysis, understanding the mathematical nuances of set representations is paramount. This section presents a concise exploration into three pivotal set representations: polyhedral sets, zonotopes, and the more intricate constrained zonotopes.

\noindent \textbf{Polyhedral Set:} Polyhedral sets, often encountered in linear programming and optimization, emerge as the intersection of a limited number of halfspaces. For matrices $\mathbf{H} \in \mathbb{R}^{m \times n}$ and vectors $\mathbf{h} \in \mathbb{R}^{m}$, the polyhedral set in $\mathbb{R}^{n}$, symbolized as $\mathcal{H}$, is formally articulated as:
\begin{equation}
	\mathcal{H} = { \mathbf{x} \in \mathbb{R}^{n} : \mathbf{H} \mathbf{x} \leq \mathbf{h} }.
\end{equation}

\noindent \textbf{Zonotope:} Given a center $\mathbf{c} \in \mathbb{R}^{n}$ and a generator matrix $\mathbf{G} \in \mathbb{R}^{n \times p}$, a zonotope $\mathcal{Z} \subset \mathbb{R}^{n}$ is defined as:
\begin{equation}
	\mathcal{Z} = { \mathbf{c} + \mathbf{G} \xi : |\xi|_{\infty} \leq 1 }.
\end{equation}

\noindent The definition of the constrained zonotope introduced in Reference~[\citenum{scott2016constrained}] is:

\noindent \textbf{Constrained Zonotope:} Advancing the foundational concept of zonotopes, constrained zonotopes offer a more nuanced representation by intertwining additional constraints. As elucidated in the seminal work of Reference~[\citenum{scott2016constrained}], given the parameters: a center $\mathbf{c} \in \mathbb{R}^{n}$, a generator matrix $\mathbf{G} \in \mathbb{R}^{n \times p}$, a constraint matrix $\mathbf{A} \in \mathbb{R}^{m \times p}$, and a constraint vector $\mathbf{b} \in \mathbb{R}^{m}$, the constrained zonotope is depicted as:
\begin{equation}
	\mathcal{Z} = { \mathbf{c} + \mathbf{G} \xi : |\xi|_{\infty} \leq 1, \mathbf{A} \xi = \mathbf{b} }.
\end{equation}
This structured notation allows for the succinct representation $\mathcal{Z} = (\mathbf{G}, \mathbf{c}, \mathbf{A}, \mathbf{b})$.

\noindent \textit{Note:} Constrained zonotopes display an innate property of closure under various operations such as intersection and Minkowski sum (expressed as $\oplus$), making them versatile tools for depicting diverse convex polytopes.

\section{Computation of the Reachable Set for a Nonlinear System}\label{sec:nonlinReach}
This section considers general nonlinear systems with uncertain parameters that are constrained within certain limits. The reachability of these nonlinear systems in CORA is determined by abstracting the state-space~\cite{althoff2010reachability}. A local linearization of the nonlinear system is performed using a Taylor series expansion. To simplify notation, the state and input vectors are combined into a new vector $\mathbf{z}{k}^{T}=\left[\mathbf{x}{k}^{T}, \mathbf{u}_{k}^{T}\right]$. The Taylor series of the nonlinear system Eq.~(\ref{eq:16}) can then be written as:
\begin{equation} \label{eq:taylor}
	\begin{aligned}
		{\left(\mathbf{x}_{k+1}\right)}_{i} = \;\mathbf{g}_{i}(\mathbf{x}_k, \mathbf{u}_k, \boldsymbol{w}_k) = & \; \mathbf{g}_{i}\left(\mathbf{z}_{k}^{*}, \boldsymbol{w}_k \right)+\left. \nabla_{\mathbf{z}_{k}}\biggl(\mathbf{g}_{i} \left(\mathbf{z}_{k}^{*}, \boldsymbol{w}_k \right) \biggr) \right|_{\mathbf{z}_{k}=\mathbf{z}_{k}^{*}}\left(\mathbf{z}_{k}-\mathbf{z}_{k}^{*}\right)+\\
		&\left.\frac{1}{2}\left(\mathbf{z}_{k}-\mathbf{z}_{k}^{*}\right)^{T} \nabla^2_{\mathbf{z}_{k}}\biggl(\mathbf{g}_{i} \left(\mathbf{z}_{k}^{*}, \boldsymbol{w}_k \right) \biggr)\right|_{\mathbf{z}_{k}=\mathbf{z}_{k}^{*}}\left(\mathbf{z}_{k}-\mathbf{z}_{k}^{*}\right)+\ldots
	\end{aligned}
\end{equation}
where $i$ refers to the $i^{th}$ coordinate of $g$ and $k$ is the time instant.\\
The infinite Taylor series can be over-approximated by a using its first-order representation paired with the Lagrange remainder:
\begin{equation} \label{eq:rem}
	\begin{aligned}
		{\left(\mathbf{x}_{k+1}\right)}_{i} & \in \underbrace{\mathbf{g}_{i}\left(\mathbf{z}_{k}^{*}, \boldsymbol{w}_k \right)+\left.\nabla_{\mathbf{z}_{k}}\biggl(\mathbf{g}_{i} \left(\mathbf{z}_{k}^{*}, \boldsymbol{w}_k \right) \biggr) \right|_{\mathbf{z}_{k}=\mathbf{z}_{k}^{*}}\left(\mathbf{z}_{k}-\mathbf{z}_{k}^{*}\right)}_{1^{s t} \text { order Taylor series }}+\\
		& \;\;\;\;\underbrace{\frac{1}{2}\left(\mathbf{z}_{k}-\mathbf{z}_{k}^{*}\right)^{T} \left.\nabla^2_{\mathbf{z}_{k}}\biggl(\mathbf{g}_{i} \left(\mathbf{z}_{k}^{*}, \boldsymbol{w}_k \right) \biggr) \right|_{\pmb{\zeta}}\left( \mathbf{z}_{k}-\mathbf{z}_{k}^{*}\right)}_{\text {Lagrange remainder } (L_{k})_{i}}
	\end{aligned}
\end{equation}
where $\pmb{\zeta} \in\left\{\mathbf{z}_{k}^{*}+\alpha\left(\mathbf{z}_{k}-\mathbf{z}_{k}^{*}\right) \mid \alpha \in[0,1]\right\}$.

The reachable set for the linearized system ($\mathcal{R}^{lin}(t)$) is first computed without accounting for linearization errors. The reachable set due to linearization errors ($\mathcal{R}^{err}(t)$) is then computed by evaluating the Lagrange remainder (Eq.~(\ref{eq:rem})). The overall reachable set ($\mathcal{R}(t)$) is obtained by performing Minkowski addition of the two sets, given by
\begin{equation}
	\mathcal{R}(t) = \mathcal{R}^{lin}(t) \oplus \mathcal{R}^{err}(t)
\end{equation}
\subsection{Computation of the Set of Linearization Errors}
\begin{equation}
	\mathbf{J}_{i}(\pmb{\zeta}, \boldsymbol{w}_k) = \left.\nabla^2_{\mathbf{z}_{k}}\biggl(\mathbf{g}_{i} \left(\mathbf{z}_{k}^{*}, \boldsymbol{w}_k \right) \biggr) \right|_{\pmb{\zeta}}
\end{equation}

where $i$ refers to the $i^{th}$ coordinate of $\mathbf{g}$, we can write the Lagrange remainder as 
\begin{equation}
	\begin{aligned}
		&\left(L_{k}\right)_{i}=\frac{1}{2}\left(\mathbf{z}_{k}-\mathbf{z}_{k}^{*}\right)^{T} \mathbf{J}_{i}(\pmb{\zeta}, \boldsymbol{w}_k)\left(\mathbf{z}_{k}-\mathbf{z}_{k}^{*}\right) \\
		&\;\;\pmb{\zeta} \in\left\{\mathbf{z}_{k}^{*}+\alpha\left(\mathbf{z}_{k}-\mathbf{z}_{k}^{*}\right) \mid \alpha \in[0,1]\right\}
	\end{aligned}
\end{equation}

The absolute values of the Lagrange remainder can be over-approximated for $\mathbf{z}_{k} \in {\mathcal{Z}}$ where ${\mathcal{Z}}$ is a zonotope ${\mathcal{Z}}=\left(\mathbf{c}, \mathbf{g}^{(1)}, \ldots, \mathbf{g}^{(e)}\right)$. The over-approximation is obtained by the following computations:
\begin{equation}
	\left|\left(L_{k}\right)_{i}\right| \subseteq\left[0,\left( l_{k}\right)_{i}\right]
\end{equation}

with $\left( l_{k}\right)_{i}=\frac{1}{2} \gamma^{T} \max \left(\left|\mathbf{J}_{i}\left(\pmb{\zeta}\left(\mathbf{z}_{k}, \mathbf{z}_{k}^{*}\right), \boldsymbol{w}_k \right)\right|\right) \gamma,  \mathbf{z} \in \mathcal{Z}$ and $\gamma=\left|\mathbf{c}-\mathbf{z}_{k}^{*}\right|+\sum_{i}^{e}\left|\mathbf{g}^{(i)}\right|$, where $\mathbf{c}$ is the center and $\mathbf{g}^{(i)}$ are the generators of the zonotope $\mathcal{Z}$.

\noindent The reachable set is computed iteratively for smaller time intervals $t \in[(k-1) \Delta t, k \Delta t]$ where $k \in$ $\mathbb{N}^{+}$, such that
$\mathcal{R}\left(\left[0, t_{f}\right]\right)$ is obtained by their union: 
\begin{equation}
	\mathcal{R}\left(\left[0, t_{f}\right]\right)=\bigcup\limits_{k=1}^{t_{f} / \Delta t} \mathcal{R}([(k-1) \Delta t, k \Delta t])
\end{equation}
At every time step $k$, the nonlinear system undergoes a linear approximation through Taylor series expansion. The reachable set for this linear model is then computed. To maintain precision, any linearization errors are restricted to a user-specified tolerance. Should these errors surpass the predefined threshold, the reachable set is bifurcated into two segments. Computations are then reiterated for each segment. This division is executed along a hyperplane, oriented perpendicular to the axis which most significantly reduces the linearization error.

\section{Heat Transfer Analysis}
During re-entry, vehicles undergo intense heating rates that are highly dependent on burnout conditions. If the burnout angle is steep, the rate of descent will be rapid, causing the vehicle to quickly reach the heat rate boundary. This can make the trajectory unfeasible due to exceeding the heat rate limit. The heating rate at the stagnation point can be calculated using the equation presented in Reference~{[\citenum{Scott1984}]}. The convective heat transfer rate is given by the equation:
\begin{equation}\label{heat_rate}
	\dot{Q}  =  C\; \rho^{0.5} v^{3.05}
\end{equation}
where the heating rate, $\dot{Q}$, is directly proportional to the density $\rho^{0.5}$ and velocity $v^{3.05}$. The proportionality constant $C$ in the equation for the stagnation point heating rate depends on the size of the re-entry vehicle's nose radius ($R_{N}$). Furthermore, the total heat load at the stagnation point is calculated using the equation:
\begin{equation}
	Q = \int_{t_{0}}^{t_{f}} \dot{Q} \, dt
\end{equation}
where $t_{0}$ and $t_{f}$ are the initial and final times, respectively.

\subsection{Heat Rate Constraint}
\begin{figure}[h]
  \centering
 \includegraphics[trim=0 0 0 0, width=0.8\columnwidth]{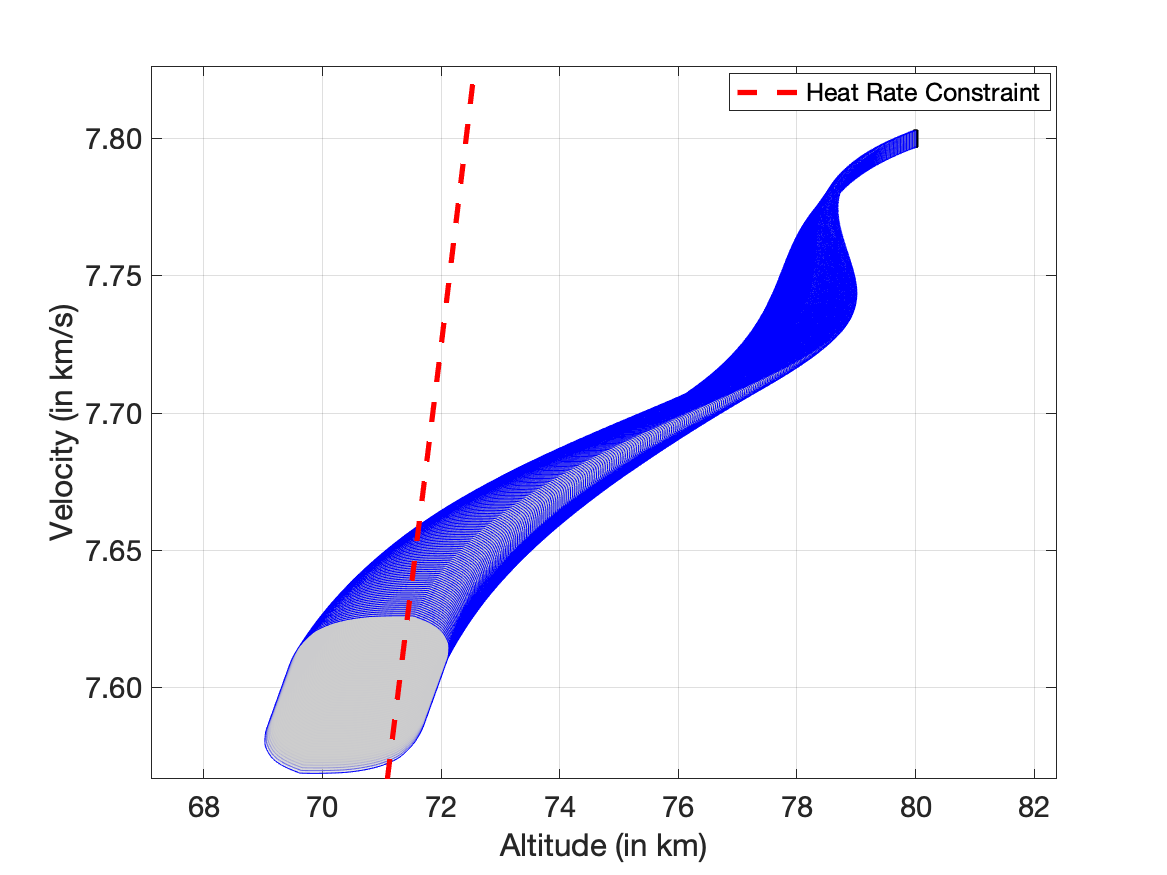}  
  \caption{Reachable set tube (Altitude vs Velocity) and heat rate limit }
  \label{fig:heat}
\end{figure}

The surface temperature of a re-entry vehicle is closely related to its heat rate limit, assuming that the surface is in thermal equilibrium with its surroundings through radiation. 
\begin{equation}
	\dot{Q} = \sigma \; \epsilon \; (T_{s}^{4})
\end{equation}
\begin{equation}
	T_{s} = \left(\frac{\dot{Q}}{\sigma \; \epsilon}\right)^{(1/4)}
\end{equation}
where, $\sigma$ is the Stefan – Boltzmann constant, which is equal to
$5.67 \times  10^{-8} \; W/(m^2.K^4)$, and $\epsilon$ is the surface emissivity which is assumed to be 0.8. To ensure that the re-entry vehicle does not exceed its heat rate limit, the velocity-altitude reachable set in Fig.~\ref{fig:heat} takes into account the constraint on the heating rate. The value of the constant $C$ in the equation Eq.~(\ref{heat_rate}) is determined to be $1.1813 \times 10^{-3} \; kg^{0.5} m^{1.5}/s^{3}$. The re-entry vehicle is made of a Carbon-Carbon composite material, which can withstand temperatures up to $2800^{\circ} C$,as detailed in Reference~{[\citenum{ulIslamRizvi2015}]}. This temperature equates to a heat rate limit of around $4.0 \; MW/m^{2}$, which is depicted by the dashed red line in the figure. Through this analysis, it is possible to determine if the re-entry vehicle has reached its heating limit. The reachable set's tube in the Fig.~\ref{fig:heat} contains all possible trajectories that the vehicle can take, given various control inputs within a certain range. Therefore, if the vehicle approaches its heating limit at any point during re-entry, it is possible to redirect its path within the tube to avoid exceeding the heat rate limit. In Fig.~\ref{fig:heat}, the initial state is considered bounded as a Zonotope with center given by,
\[
	\mathbf{x}=[80,000 (m), 7800 (m/s), -0.1^{\circ}, 0^{\circ}, 90^{\circ}, 0^{\circ}]
\]
The control terms $\alpha$ and $\beta$ are assumed to be bounded in intervals, $\alpha \in [15^{\circ},30^{\circ}]$ and $\beta \in [-65^{\circ},-55^{\circ}]$.

\section{Model Predictive Control}\label{mpc}
The essence of Model Predictive Control (MPC) lies in numerically optimizing a model to derive a control input sequence. This sequence aims to minimize a cost or objective function over a predefined control horizon. At each sampling instant, these control actions are updated based on the current state estimate, ensuring an adaptive feedback mechanism.
We define the state vector as \(\mathbf{x} = [h \; v \; \gamma \; \theta \; \psi \; \phi]^{\top}\) and the associated control actions as \(\mathbf{u} = [\alpha, \beta ]^{\top}\). The system dynamics for the reentry problem are discretized with a sampling time \(T_s\) through the forward Euler method, resulting in the relation \(\mathbf{x}_{k+1} = \mathbf{g}(\mathbf{x}_k, \mathbf{u}_k)\). This discrete model is integral to the predictive component of the MPC.
Predictions within the MPC span a receding horizon, meaning they encompass a defined number of future steps. This duration is denoted as the prediction horizon, \(N\).

\subsection{Cost Function}
The primary objective of the MPC is to minimize the tracking errors between the predicted states and the desired trajectory.  A cost function $J$ is form as below.
\begin{equation} 
	\begin{array}{c}
		J(k) =  \sum_{t=k}^{k+N-1} \left[ (\mathbf{x}(t) - \mathbf{x}^{\text{ref}}(t))^T \mathbf{Q} (\mathbf{x}(t) - \mathbf{x}^{\text{ref}}(t)) + (\mathbf{u}(t) - \mathbf{u}^{\text{ref}}(t))^T \mathbf{R} (\mathbf{u}(t) - \mathbf{u}^{\text{ref}}(t)) \right]
	\end{array}
\end{equation}
where $\mathbf{Q}  \in \mathbb{R}^{6\times 6}$ and $\mathbf{R} \in \mathbb{R}^{2\times 2}$ are symmetric positive definite matrices for the states and inputs, respectively. The first term denotes the state cost, which penalizes deviating from a certain state reference $\mathbf{x}^{ref}$. The second term denotes the input cost that penalizes a deviation from the steady-state input $\mathbf{u}^{ref}$. Finally, to enforce smooth control actions, bounds on the rate of control changes are also imposed.

\subsection{Constraints}
\subsubsection{Reachable Safe Sets:}
The reachable sets computed using the method from above sections are considered safe as they include all the states which are achievable given the bounds on the control inputs. In order to find a safe manuever, the re-entry vehicle should be within the reachable sets. Zonotopes have demonstrated exceptional computational efficiency but they are not closed under the intersection, thus, it is necessary to convert them to polytopes and back again, which involoves approximations. Additionally, most systems have state constraints, which should keep the reachable sets within the bounds of the state constraints. Now, if we consider the heating rate constraint that the trajectory should avoid getting past the heating rate limit. Therefore, to confirm safety specifications, the intersection of the reachable set with the constraints (which can be heating rate constraint) is performed. These state constraints can be represented as a combination of halfspacess~\cite{althoff2012avoiding}. The intersection between the zonotope and halfspace can be represented as a constrained zonotope. Constrained zonotope enables exact representations, unlike current methods that rely on zonotopic approximations of the intersection.

\subsubsection{Control Inputs:}
Some physical limitations need to be considered when computing the optimal control inputs,
\begin{equation}
	\boldsymbol{u}_k \in \mathcal{U}
\end{equation}
where $\mathcal{U}$ denotes the bounded zonotopic set.

\subsection{Implementation}
To achieve the desired calculations for the re-entry vehicle's reachable safe states, we utilized a combination of Matlab toolboxes. The fundamental mathematics for these calculations primarily hinge on set operations. Following toolboxes are utilized:

\subsubsection{MATLAB's MPC Toolbox:} Essential for incorporating constraints within the Model Predictive Control (MPC) framework. In MPC, a cost function is minimized over a finite control horizon while simultaneously taking constraints into account. 
\subsubsection{MPT3 Toolbox~\cite{MPT3}:} This toolbox allows for the representation of polyhedral sets derived from zonotopic reachable sets. These polyhedra integrate seamlessly with predefined set operations within the toolbox.
\subsubsection{YALMIP Toolbox~\cite{lofberg_yalmip_2004}:} Used for diverse optimization and modeling tasks.
\subsubsection{CORA Toolbox~\cite{althoff_implementation_2016}:} We utilized certain components of this toolbox in our implementation.\\\\
Additionally, we employ Algorithms~\ref{alg_one} and \ref{alg_two} in our approach. Algorithm~\ref{alg_one} is instrumental in computing the reachable sets, providing a comprehensive mapping of all feasible states the system can attain. On the other hand, Algorithm~\ref{alg_two} is tailored for the application of Model Predictive Control (MPC), allowing for meticulous trajectory determination that adheres to the system's constraints and desired objectives.

\begin{algorithm}
	\caption{Computation of Reachable Set with State Bounds}\label{alg_one}
	\KwData{Constrained zonotope $\mathcal{R}(0) = \mathcal{X}_{0}$, time step size $\Delta t$, and time horizon $t_f$}
	
	\KwResult{Reachable set up to time horizon $t_f$, i.e. 	$\mathcal{R}\left(\left[0, t_{f}\right]\right)$}
	
	\textbf{Initialize:} $k \gets 1$\;
	
	\While{$k \Delta t \leq t_f$}{
		Compute the reachable set for the current time step using Algorithm 1~\cite{jinay}\;
		
		$\mathcal{R}(k\Delta t) = \mathcal{H}_{-} \cap \mathcal{R}(k\Delta t)$\;
		
		\texttt{reduceorder}$(\mathcal{R}(k\Delta t))$ \tcp*[r]{Apply order reduction, see~\cite{scott2016constrained}}
		
		Update initial set to $\mathcal{R}(k \Delta t)$\;
		
		$k \gets k + 1$\;
	}
\end{algorithm}

\begin{algorithm}
	\caption{Nonlinear Model Predictive Control with Reachable Set Constraints}\label{alg_two}
	\KwData{Nonlinear system dynamics \(\mathbf{g}(\mathbf{x},\mathbf{u})\), initial state \(\mathbf{x}_0\), reference trajectories \(\{\mathbf{x}^{\text{ref}}, \mathbf{u}^{\text{ref}}\}\), reachable set $\mathcal{R}\left(\left[0, t_{f}\right]\right)$, control input bound, prediction horizon \(N\).}
	\KwResult{Optimal control sequence \(\{\mathbf{u}^*(t)\}_{t=0}^{N-1}\)}
	
	\textbf{Initialize:} \(k \gets 0\), \(\mathbf{x}(k) = \mathbf{x}_0\)\;
	
	\While{End criteria not met}{
		\textbf{Optimization problem:}
		\begin{equation*}
			\begin{aligned}
				& \underset{\{u(t)\}_{t=k}^{k+N-1}}{\text{minimize}}
				& & \sum_{t=k}^{k+N-1} \left[ (\mathbf{x}(t) - \mathbf{x}^{\text{ref}}(t))^T \mathbf{Q} (\mathbf{x}(t) - \mathbf{x}^{\text{ref}}(t)) + (\mathbf{u}(t) - \mathbf{u}^{\text{ref}}(t))^T \mathbf{R} (\mathbf{u}(t) - \mathbf{u}^{\text{ref}}(t)) \right] \\
				& \text{subject to}
				& & \mathbf{x}(t+1) = \mathbf{g}(\mathbf{x(t)},\mathbf{u(t)}), \; t = k, \ldots, k+N-1 \\
				& & & \mathbf{x}(t) \in \mathcal{R}\left(\left[0, t_{f}\right]\right) \;\;\;\;\;\;\;\;\; // \;\; \mathcal{R}\left(\left[0, t_{f}\right]\right) \text{ is computed from Algorithm 1} \\ 
				& & & \mathbf{u}(t) \in \mathcal{U}
			\end{aligned}
		\end{equation*}
		
		Solve the nonlinear optimization problem to obtain the optimal control sequence \(\{\mathbf{u}^*(t)\}_{t=k}^{k+N-1}\)\;
		
		Apply the first control action \(\mathbf{u}^*(k)\) to the system\;

		Increment time step: \(k \gets k + 1\)\;
	}
	
\end{algorithm}

\section{Simulation Results}
Numerical simulations are carried out to compute reachable sets for the reentry problem under consideration. All the simulations are performed on a computer with 16.00 GB RAM and 2.10 GHz Intel Core i7 processor running MATLAB R2022a. We present simulation results pertaining to the reentry problem under investigation, showcasing the efficacy of our model and control strategy. A summary of the parameters and their corresponding values employed in the simulation can be found in Table \ref{table:parameters}. 
\begin{table}[h]
	\centering
	\caption{Parameters used in the simulation and aerodynamic models.}
	\begin{tabular}{ccc}
		\toprule\hline
		Parameter & Value & Unit \\
		\midrule
		$G$ & $3.986 \times 10^{14}$ & m$^3$/s$^2$ \\
		$\rho_0$ & 1.225 & kg/m$^3$ \\
		$R_e$ & 6378 & km \\
		$h_r$ & 7500 & m \\
		$S$ & 0.30 & m$^2$ \\
		$m$ & 340 & kg \\
		$a_0$ & -0.20704 & -- \\
		$a_1$ & 0.029244 & -- \\
		$b_0$ & 0.07854 & -- \\
		$b_1$ & -0.61592 $\times$ $10^{-2}$ & -- \\
		$b_2$ & 0.621408 $\times$ $10^{-3}$ & -- \\
		\bottomrule
	\end{tabular}
	\label{table:parameters}
\end{table}

\begin{table}[h]
	\centering
	\caption{Control terms assumed to be bounded in given intervals.}
	\begin{tabular}{cc}
		\toprule
		\hline
		Control Term & Interval (degrees) \\
		\midrule
		$\alpha$ (angle of attack) & $[15, 30]$ \\
		$\beta$ (bank angle) & $[-60, -50]$ \\
		\bottomrule
	\end{tabular}
	\label{table:control_terms}
\end{table}

\subsection{Implementing MPC}
\begin{figure}[h]
	\centering
	\includegraphics[trim=0 0 0 0, width=0.9\columnwidth]{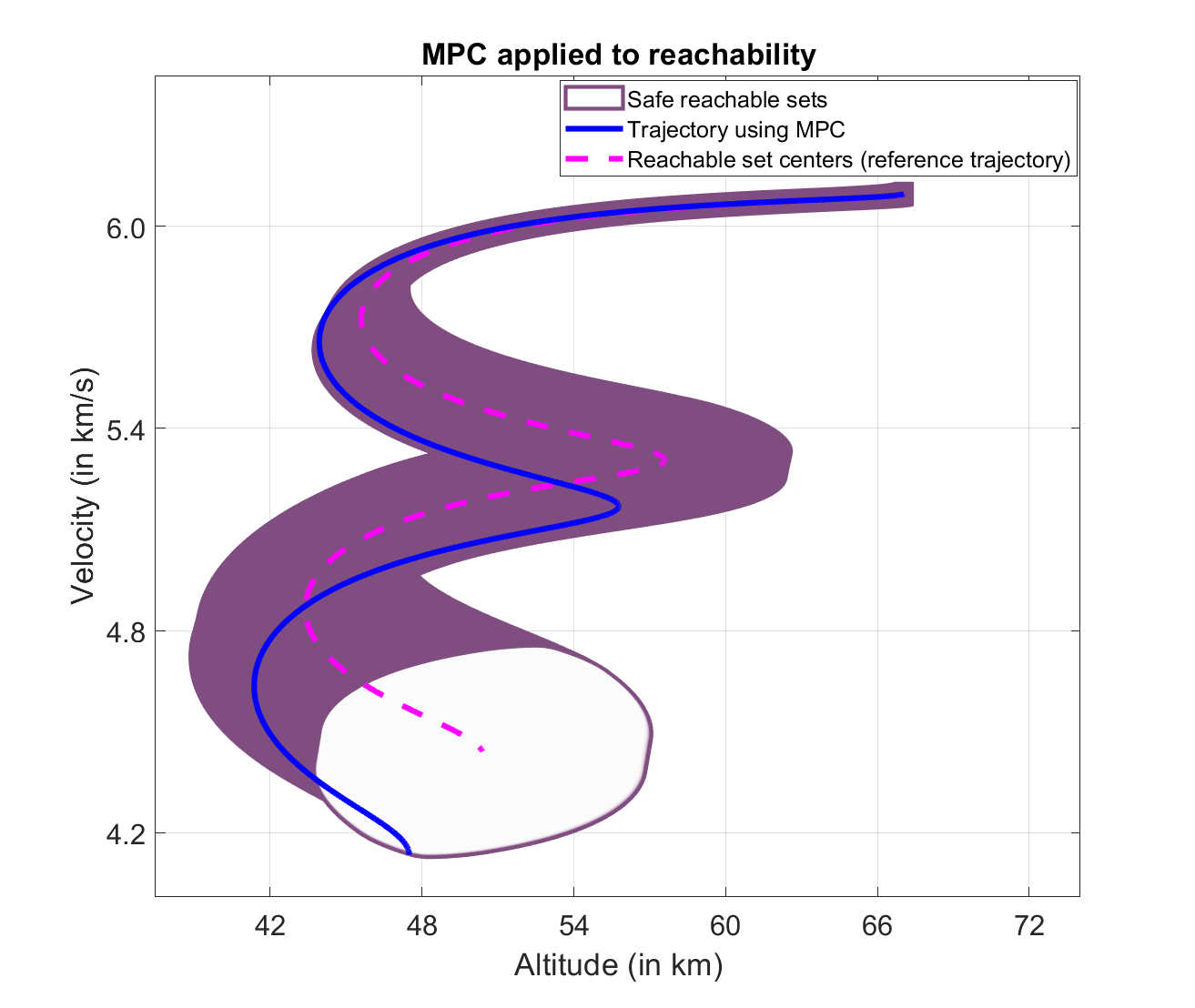}  
	\caption{Zonotope reachable sets evolution}
	\label{fig:reachset_without_heat_limit}
\end{figure}

The evolution of reachable sets is illustrated in Fig. \ref{fig:reachset_without_heat_limit}. These sets were computed based on the assumption that the control input values are bounded as outlined in Table \ref{table:control_terms}. Additionally, we determined an optimal trajectory, assuming the initial condition is within the initial reachable set, utilizing the Model Predictive Control approach detailed in Algorithm \ref{alg_two}. It's important to note that these computed reachable sets are integrated as constraints within the MPC. In Fig. \ref{fig:reachset_without_heat_limit}, the trajectory ascertained through the MPC is displayed. For our purposes, the reference trajectory was set to align with the center of the zonotopes, ensuring that we remain within the bounds of the reachable set while maintaining a satisfactory safety margin. The chosen time horizon, $t_f$, stands at $1200~\mathrm{s}$ and the sampling time is set at $0.1~\mathrm{s}$. The initial state is defined as:
\[
\mathbf{x}=[67,000 (m), 6096 (m/s), -0.1^{\circ}, 0^{\circ}, 90^{\circ}, 0^{\circ}]
\]
Furthermore, the MPC prediction horizon $N$ has been configured with $N=20$, and the weights are determined as $\mathbf{Q} = 100 \; \mathbf{I}_{6\times 6}$ and $\mathbf{R} =  \mathbf{I}_{2\times 2}$.

\subsection{Addressing Heat Rate Constraint}
\begin{figure}[h]
  \centering
 \includegraphics[trim=0 0 0 0, width=0.8\columnwidth]{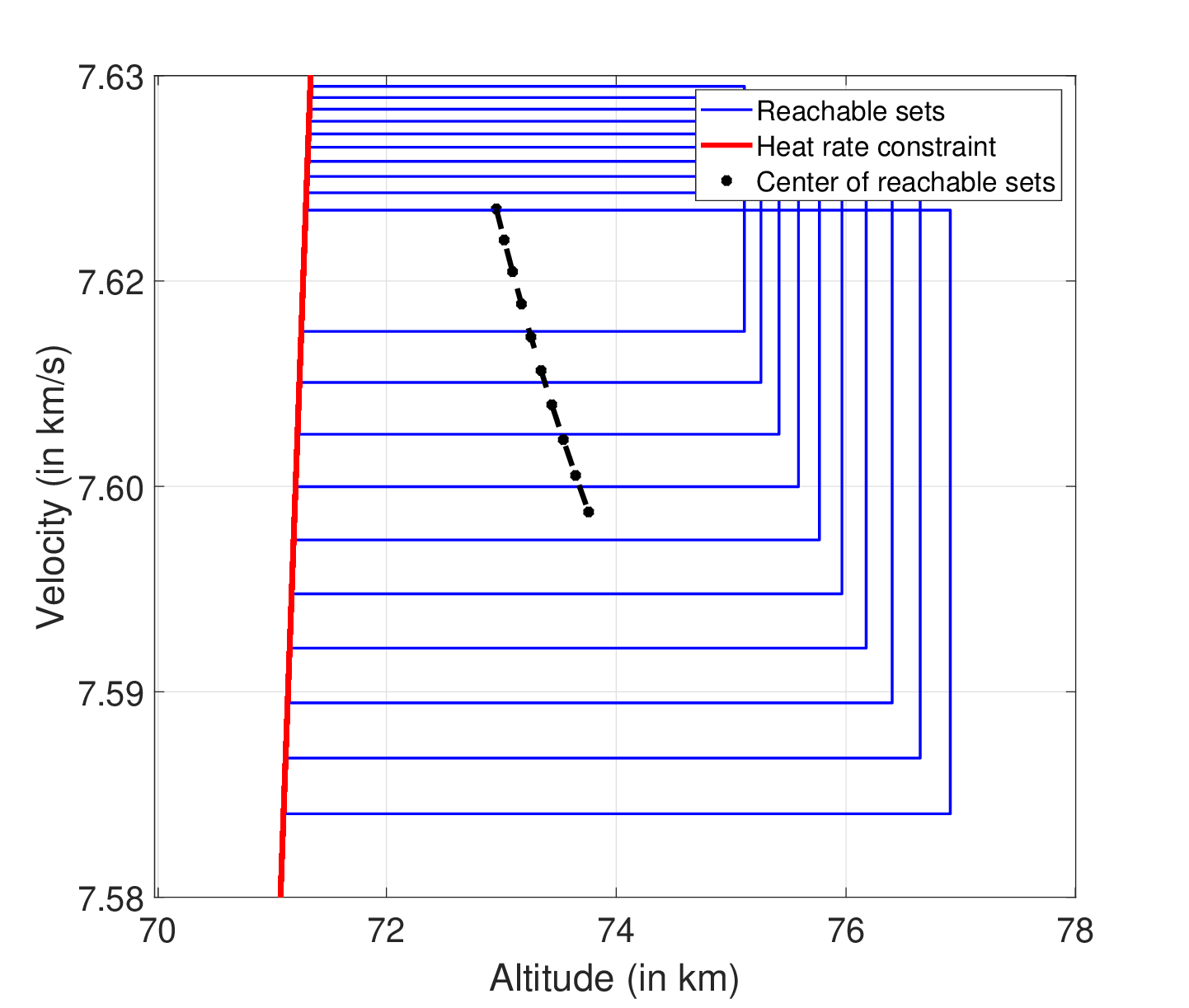}  
  \caption{Constrained zonotope reachable sets evolution given the heating limit (red-line) }
  \label{fig:conZono}
\end{figure}

To handle the heat rate constraint illustrated in Fig.~\ref{fig:heat}, we will employ constrained zonotopes, as proposed in Lemma 2~\cite{jinay}. This involves formulating the heat rate constraint as a hyperplane and intersecting it with a zonotope to yield a constrained zonotope. Fig.~\ref{fig:conZono} showcases the use of constrained zonotopes in implementing the heat constraint. To begin with, an initial state is defined as a Zonotope with a center given by
\[
	\mathbf{x}=[71,932 (m), 7600 (m/s), -0.1^{\circ}, 0^{\circ}, 90^{\circ}, 0^{\circ}]
\]

Furthermore, we assume that the control terms $\alpha$ and $\beta$ are bounded in intervals, $\alpha \in [15^{\circ},30^{\circ}]$ and $\beta \in [-60^{\circ},-10^{\circ}]$.

\begin{figure}
	\centering
	\begin{minipage}{0.49\textwidth}
		\includegraphics[trim=0 0 0 0, width=1\columnwidth]{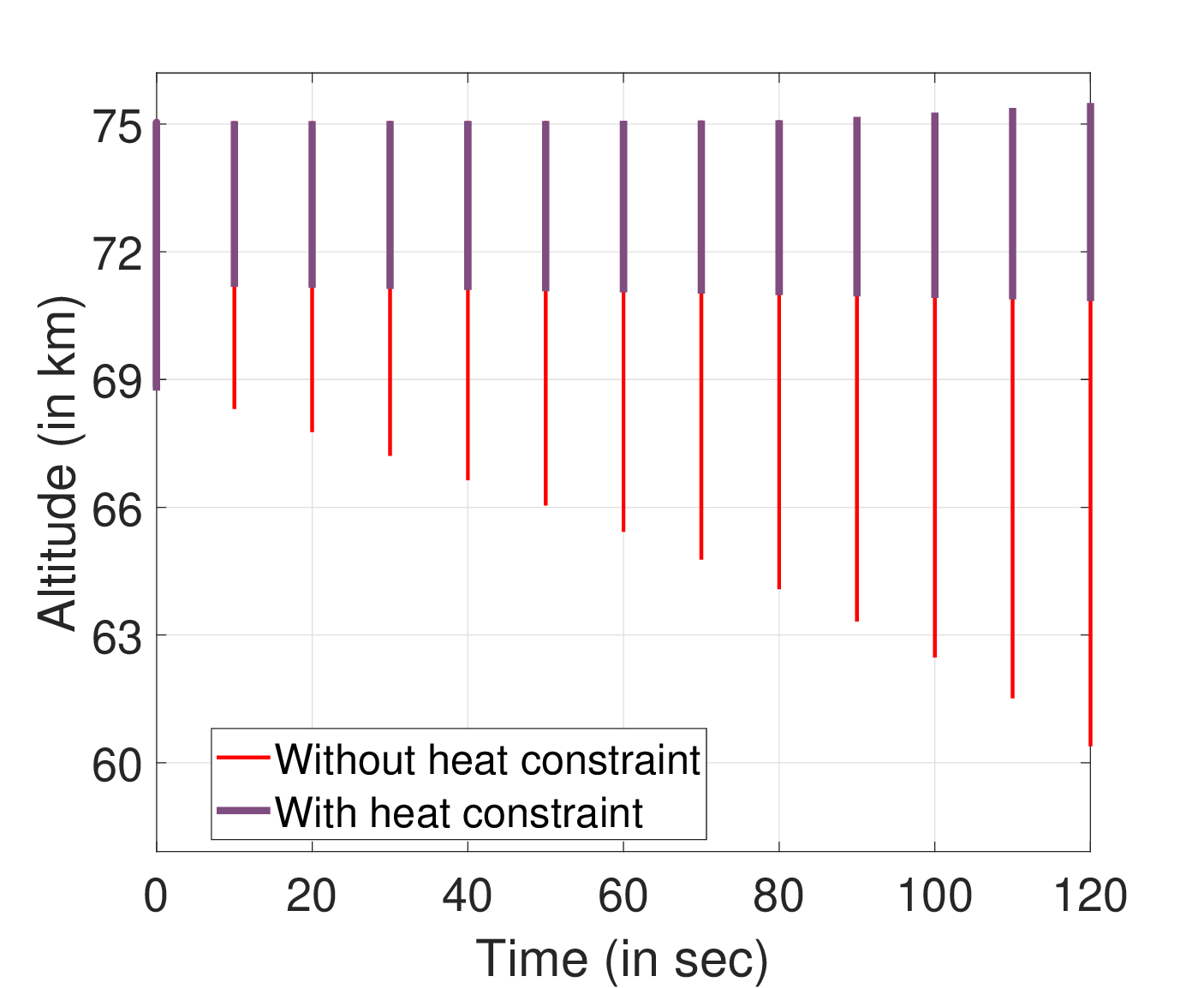}  
		\subcaption{Altitude vs Time }
		\label{fig:alt}
	\end{minipage}
	\begin{minipage}{0.49\textwidth}
		\includegraphics[trim=0 5 0 0, width=1\columnwidth]{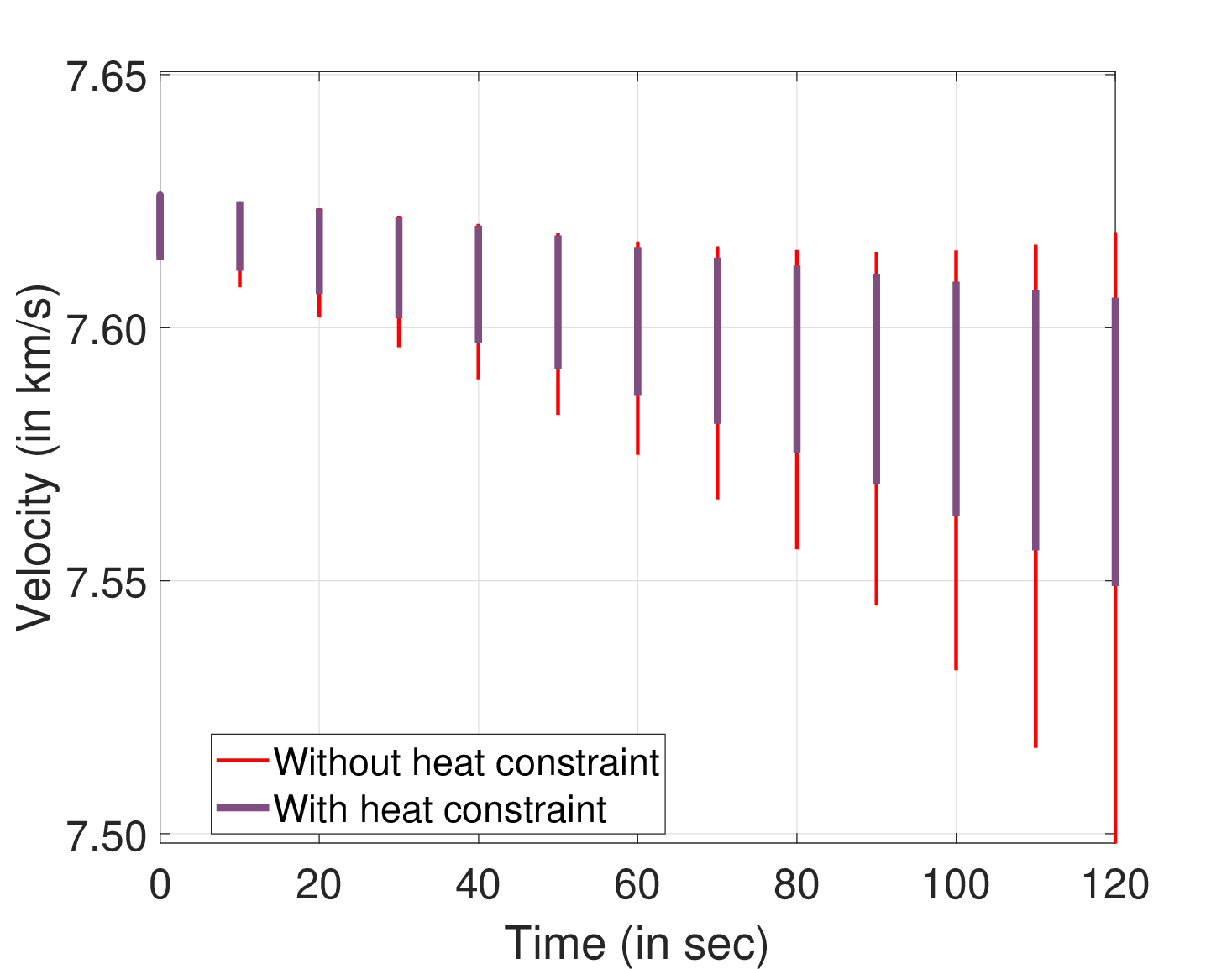}  
		\subcaption{Velocity vs Time }
		\label{fig:vel}
	\end{minipage}
	\caption{Constrained zonotope evolution with time}
	\label{fig:alt-vel}
\end{figure}

Figure~\ref{fig:alt-vel} illustrates the evolution of the constrained zonotope with time for two cases: without any heating limit constraint and with a heating limit constraint. It is observed that when the constraint is imposed on the heating load, the altitude does not decrease, as the heat rate increases as the altitude drops. As a result, the reachable sets for the altitude dimension remain the same, as shown in Figure~\ref{fig:alt}. However, for velocity, the reachable sets start to decrease, which reduces the heating load and moves the vehicle away from the heating limit. The reachable sets were obtained using a time step size of $\Delta t = 10~\mathrm{s}$ and a time horizon of $t_f = 120~\mathrm{s}$.

\subsection{Heat-Constrained Trajectory Planning}
\begin{figure}[h]
	\centering
	\includegraphics[trim=0 0 0 0, width=0.8\columnwidth]{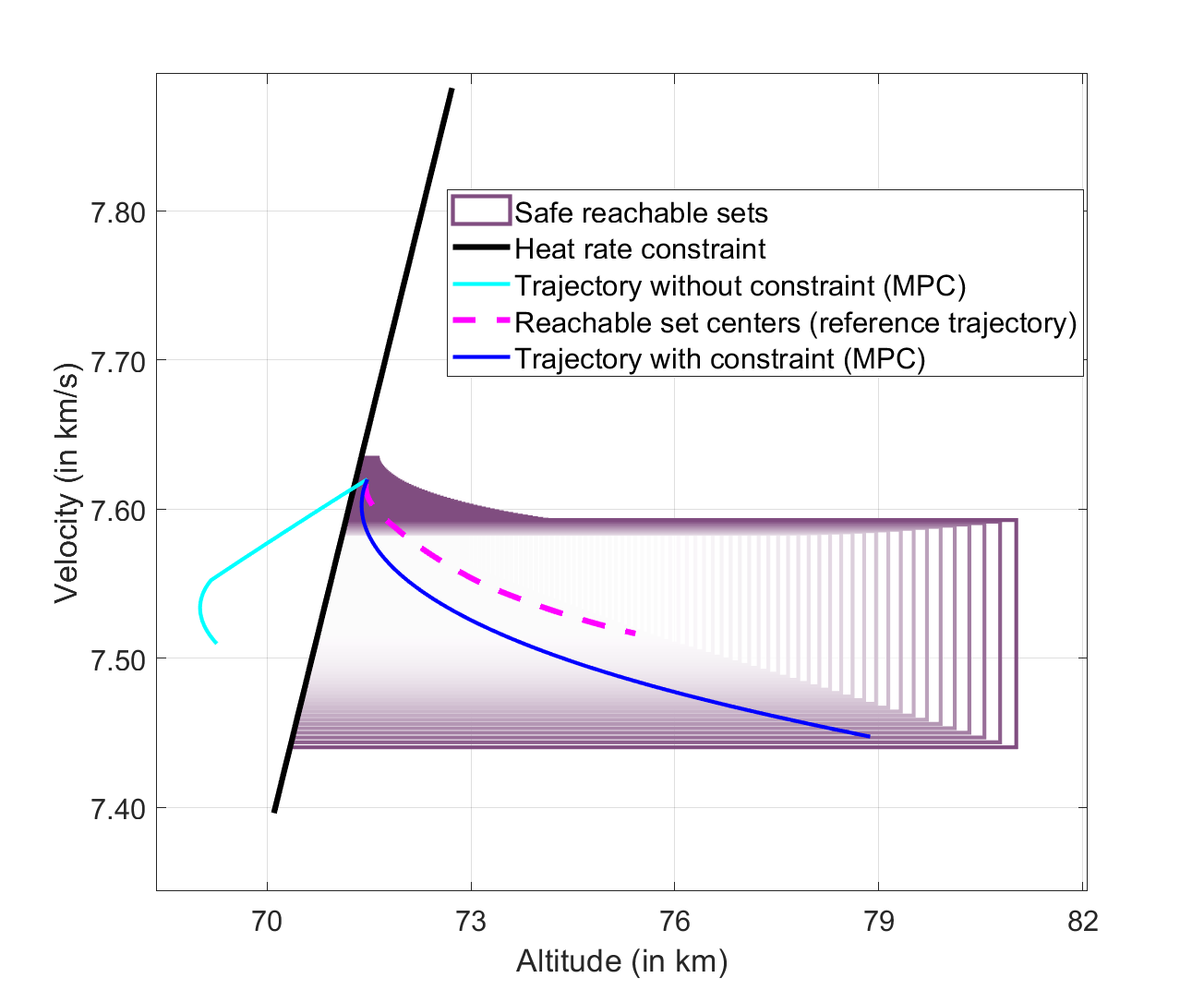}  
	\caption{Constrained zonotope reachable sets evolution given the heating limit (red-line) }
	\label{fig:conZono2}
\end{figure}

This subsection delves into the trajectory planning, particularly addressing the heat rate constraint. The starting zonotope is defined as:
\[
\mathbf{x}=[71,171 (m), 7620 (m/s), -0.1^{\circ}, 0^{\circ}, 90^{\circ}, 0^{\circ}]
\]
Reachable sets, based on the assumption of bounded control input values, are derived as detailed in Table \ref{table:control_terms2}. The procedure leverages Algorithm~\ref{alg_one} to implement the heat rate constraint. The representation of the heat rate constraint employs a halfspace, wherein the intersection between a zonotope and a halfspace yields a constrained zonotope. This constrained zonotope is then evolved. In subsequent iterations, the new intersection between this evolved constrained zonotope and a halfspace produces another constrained zonotope. The simulations utilize a time step of $0.5~\mathrm{s}$. In Fig.~\ref{fig:conZono2}, the evolution of altitude-velocity reachable sets over time is vividly captured. Among the trajectories presented, two distinct paths can be discerned. The trajectory delineated by the cyan line represents a scenario where the heat rate limit is disregarded. On the other hand, the trajectory highlighted in blue takes the heat rate constraint into consideration, notably maintaining its course within the initial zonotope prior to intersecting with the halfspace. This visual representation underscores the precision of the Model Predictive Control approach: while it strategically navigates around the heat rate limit, it simultaneously ensures the trajectory conforms to the predefined reachable sets, serving as inherent constraints in the MPC design.

\begin{table}[h]
	\centering
	\caption{Control terms assumed to be bounded in given intervals.}
	\begin{tabular}{cc}
		\toprule
		\hline
		Control Term & Interval (degrees) \\
		\midrule
		$\alpha$ (angle of attack) & $[15, 30]$ \\
		$\beta$ (bank angle) & $[-60, -30]$ \\
		\bottomrule
	\end{tabular}
	\label{table:control_terms2}
\end{table}

\section{Summary and Conclusion}

This paper introduces reachability analysis for the re-entry problem and presents preliminary results. The objective is to characterize the reachable sets of states and outputs using control variables such as the angle of attack and bank angle. To improve the practicality of set-based techniques, the paper employs zonotopes and constrained zonotopes for set operations, which offer significant computational advantages. It has been demonstrated that a constrained zonotope framework is especially useful when state constraints are present, as it provides more accurate results due to being closed under an intersection. In addition, a heat analysis is performed to prevent the vehicle from exceeding the heat constraint. The Model Predictive Control is implement to generate the re-entry vehicle trajectory. In the event that the re-entry vehicle approaches its maximum allowable heat level, MPC is employed to guide it back into the reachability tube and maintain its position within the tube. This will provide confidence that any reference trajectory can be achieved while maintaining a safe heating rate.

\bibliographystyle{AAS_publication}   
\bibliography{references}   

\begin{thebibliography}{10}

\bibitem{arslantacs2016safe}
Y.~E. Arslanta{\c{s}}, T.~Oehlschl{\"a}gel, and M.~Sagliano, ``Safe landing
  area determination for a moon lander by reachability analysis,''  {\em Acta
  Astronautica}, Vol.~128, 2016, pp.~607--615.

\bibitem{rajnishReentry}
R.~Bhusal and K.~Subbarao, ``Nonlinear State Estimation of Re-entry Vehicle
  Using Polynomial Chaos-Based Ensemble Filtering,''  {\em AAS/AIAA
  Astrodynamics Specialists Conference}, No.~AAS-20-578.

\bibitem{katiReentry}
K.~Balachandran and K.~Subbarao, ``Optimal Trajectory Control of a Space
  Shuttle like Re-entry Vehicle entering Triton,''  {\em AAS/AIAA Astrodynamics
  Specialists Conference}, No.~AAS 20-771.

\bibitem{zhou2015reachable}
Y.~Zhou and J.~S. Baras, ``Reachable set approach to collision avoidance for
  uavs,''  {\em 2015 54th IEEE Conference on Decision and Control (CDC)}, IEEE,
  2015, pp.~5947--5952.

\bibitem{malone2017hybrid}
N.~Malone, H.-T. Chiang, K.~Lesser, M.~Oishi, and L.~Tapia, ``Hybrid dynamic
  moving obstacle avoidance using a stochastic reachable set-based potential
  field,''  {\em IEEE Transactions on Robotics}, Vol.~33, No.~5, 2017,
  pp.~1124--1138.

\bibitem{mcmahon2014sampling}
T.~McMahon, S.~Thomas, and N.~M. Amato, ``Sampling-based motion planning with
  reachable volumes: Theoretical foundations,''  {\em 2014 IEEE International
  Conference on Robotics and Automation (ICRA)}, IEEE, 2014, pp.~6514--6521.

\bibitem{lofberg_yalmip_2004}
J.~Lofberg, ``{YALMIP} : a toolbox for modeling and optimization in {MATLAB},''
   {\em 2004 {IEEE} {International} {Conference} on {Robotics} and
  {Automation}}, Sept. 2004, pp.~284--289, 10.1109/CACSD.2004.1393890.

\bibitem{girard2005reachability}
A.~Girard, ``Reachability of uncertain linear systems using zonotopes,''  {\em
  International Workshop on Hybrid Systems: Computation and Control}, Springer,
  2005, pp.~291--305.

\bibitem{althoff_implementation_2016}
M.~Althoff and D.~Grebenyuk, ``Implementation of {Interval} {Arithmetic} in
  {CORA} 2016,''  Apr. 2016, 10.29007/w19b.

\bibitem{kurzhanskiy_ellipsoidal_2006}
A.~A. Kurzhanskiy and P.~Varaiya, ``Ellipsoidal toolbox (ET),''  {\em
  Proceedings of the 45th IEEE Conference on Decision and Control}, IEEE, 2006,
  pp.~1498--1503.

\bibitem{althoff_reachability_2008}
M.~Althoff, O.~Stursberg, and M.~Buss, ``Reachability analysis of nonlinear
  systems with uncertain parameters using conservative linearization,''  {\em
  2008 47th IEEE Conference on Decision and Control}, IEEE, 2008,
  pp.~4042--4048.

\bibitem{filippova_ellipsoidal_2017}
T.~F. Filippova, ``Ellipsoidal estimates of reachable sets for control systems
  with nonlinear terms,''  {\em IFAC-PapersOnLine}, Vol.~50, No.~1, 2017,
  pp.~15355--15360.

\bibitem{schurmann2018reachset}
B.~Sch{\"u}rmann, N.~Kochdumper, and M.~Althoff, ``Reachset model predictive
  control for disturbed nonlinear systems,''  {\em 2018 IEEE Conference on
  Decision and Control (CDC)}, IEEE, 2018, pp.~3463--3470.

\bibitem{7330732}
C.~A. Pascucci, S.~Bennani, and A.~Bemporad, ``Model predictive control for
  powered descent guidance and control,''  {\em 2015 European Control
  Conference (ECC)}, 2015, pp.~1388--1393, 10.1109/ECC.2015.7330732.

\bibitem{jinay_JGCD}
J.~Patel and K.~Subbarao, ``Reachability Analysis for Atmospheric Reentry
  Vehicle,''  {\em Journal of Guidance, Control, and Dynamics}, Vol.~47, No.~1,
  2024, pp.~133--142, 10.2514/1.G007549.

\bibitem{graichen2008constructive}
K.~Graichen and N.~Petit, ``Constructive methods for initialization and
  handling mixed state-input constraints in optimal control,''  {\em Journal Of
  Guidance, Control, and Dynamics}, Vol.~31, No.~5, 2008, pp.~1334--1343.

\bibitem{jain2020computationally}
A.~Jain, D.~Gueho, and P.~Singla, ``A computationally efficient approach for
  stochastic reachability set analysis,''  {\em AIAA Scitech 2020 Forum}, 2020,
  p.~0851.

\bibitem{kailath}
T.~Kailath, {\em Linear Systems}.
\newblock Prentice-Hall Inc., Englewood Cliffs, N.J, USA, 1980.

\bibitem{scott2016constrained}
J.~K. Scott, D.~M. Raimondo, G.~R. Marseglia, and R.~D. Braatz, ``Constrained
  zonotopes: A new tool for set-based estimation and fault detection,''  {\em
  Automatica}, Vol.~69, 2016, pp.~126--136, 10.1016/j.automatica.2016.02.036.

\bibitem{althoff2010reachability}
M.~Althoff, {\em Reachability analysis and its application to the safety
  assessment of autonomous cars}.
\newblock PhD thesis, Technische Universit{\"a}t M{\"u}nchen, 2010.

\bibitem{Scott1984}
C.~D. Scott, R.~C. Ried, R.~Maraia, C.~P. Li, and S.~M. Derry, ``An AOTV
  aeroheating and thermal protection study,''  {\em 19th Thermophysics
  Conference}, 1984, 10.2514/6.1984-1710.

\bibitem{ulIslamRizvi2015}
S.~T. u.~Islam~Rizvi, H.~Linshu, and X.~Dajun, ``Optimal trajectory analysis of
  hypersonic boost-glide waverider with heat load constraint,''  {\em Aircraft
  Engineering and Aerospace Technology}, Vol.~87, Jan. 2015, pp.~67--78,
  10.1108/aeat-04-2013-0079.

\bibitem{althoff2012avoiding}
M.~Althoff and B.~H. Krogh, ``Avoiding Geometric Intersection Operations in
  Reachability Analysis of Hybrid Systems,''  {\em Proceedings of the 15th ACM
  International Conference on Hybrid Systems: Computation and Control}, New
  York, NY, USA, Association for Computing Machinery, 2012, p.~45–54,
  10.1145/2185632.2185643.

\bibitem{MPT3}
M.~Herceg, M.~Kvasnica, C.~Jones, and M.~Morari, ``{Multi-Parametric Toolbox
  3.0},''  {\em Proc.~of the European Control Conference}, Z\"urich,
  Switzerland, July 17--19 2013, pp.~502--510.
\newblock \url{http://control.ee.ethz.ch/~mpt}.

\bibitem{jinay}
J.~Patel, P.~Arora, and K.~Subbarao, ``Reachable Set Computation and Analysis
  for Hypersonic Atmospheric Re-Entry Vehicles,''  {\em AAS/AIAA Astrodynamics
  Specialists Conference}, No.~AAS 22-783, Charlotte, North Carolina, AAS
  22-783, August 2022.

\end{thebibliography}

\end{document}